\documentclass[12pt]{article}
\usepackage{geometry}                
\usepackage{graphicx}
\usepackage{amssymb}					
\usepackage{amsmath}					
\usepackage{mathrsfs}
\usepackage{tikz}
\usepackage{epstopdf}
\usepackage{amsthm}
\newtheorem*{theorem*}{Theorem}
\newtheorem*{lemma*}{Lemma}
\newtheorem*{example*}{Example}
\newtheorem{theorem}{Theorem}
\newtheorem{lemma}{Lemma}

\newtheorem{proposition}{Proposition}

\def\b0{{\bf 0}}
\def\b1{{\bf 1}}

\def\n{\noindent}

\begin{document}
\title{Cubic planar bipartite graphs are dispersable \thanks{As submitted for publication on 7/20/2020; see Postscript for additional information.}
}
\author{
Paul C. Kainen\\
 \texttt{kainen@georgetown.edu}
\\
Shannon Overbay\\
\texttt{overbay@gonzaga.edu}}
\date{}                                           

\newcommand{\Addresses}{{
  \bigskip
  \footnotesize

\n
Paul C. Kainen, \textsc{Department of Mathematics and Statistics,\\ Georgetown University, 
37th and O Streets, N.W., Washington DC 20057}\\
\\
\n
Shannon Overbay, \textsc{Department of Mathematics, Gonzaga University,\\
502 E. Boone Ave., Spokane WA 99258}\\
\vspace{-0.27cm}
\n

\par\nopagebreak
}}

\maketitle

\abstract{
\n
A graph is called {\it dispersable} if it has a book embedding in which each page has maximum degree 1 and the number of pages is the maximum degree. Bernhart and Kainen conjectured every $k$-regular bipartite graph is dispersable.  Forty years later, Alam, Bekos, Gronemann, Kaufmann, and Pupyrev
have disproved this conjecture, identifying nonplanar 3- and 4-regular bipartite graphs that are not dispersable. They also proved all cubic planar bipartite 3-connected graphs are dispersable and conjectured that the connectivity condition could be relaxed. We prove that {\it every} cubic planar bipartite {\it multi}graph is dispersable.

}

\smallskip

\n
{\bf Key Phrases}: {\it book thickness; dispersable book embeddings; matching book thickness;
subhamiltonian vertex order; cubic bipartite planar multigraphs}.

\section{Introduction}

The {\it book thickness} of a graph, \cite{pck-gw, lo-fau} is now well-studied and
 has a variety of applications such as manipulation of stacks and queues (Even \& Itai \cite{even-itai} and Tarjan \cite{tarjan}), fault-tolerant computing \cite{chung-et-al}, traffic engineering \cite{pck-90},  RNA folding  \cite{haslinger-stadler}, and graph drawing and visualization, e.g., \cite{eppstein2001, eppstein2002, wood}.  There is also an extensive literature on the purely mathematical aspects of this measure of graph complexity; \cite{wiki} lists over a couple of dozen and is by no means complete.

The book thickness of a graph is the least number of colors needed for the edges such that, for some outerplane drawing of the graph, no two edges of the same color cross.
In \cite{bk79}, we introduced the additional constraint that same-color edges  {\it do not share an endpoint}, so pages must be matchings.

A {\it book embedding} of a graph is an outerplane drawing (i.e., a cyclic order on the vertices) and an edge-coloring such that edges of each color (called a {\it page} of the embedding) induce an outerplane subgraph. 
A graph or vertex-order with a 2-page book embedding is called {\it subhamiltonian} \cite{bk79}, \cite{so-thesis}.

The {\bf matching book thickness} \cite{pck-circ} of a graph $G$ is the least number of pages, $mbt(G)$, needed if the pages are matchings.
A graph is said to be {\bf dispersable} \cite{bk79} if it has a matching book embedding in $\Delta(G)$ pages, where $\Delta(G)$ denotes maximum degree; we refer to a {\bf dispersable book embedding}.

Bernhart and Kainen \cite{bk79} conjectured (in 1979):\\
$
\mbox{\it If $G$ is regular bipartite, then $G$ has a dispersable book embedding}.
$\\
This BK Conjecture was recently disproved by Alam et al. \cite{abgkp2018} who showed:\\
$
\mbox{\it There exist bipartite 3- and 4-regular graphs $G$ with $mbt(G)= 4$ or $5$} 
$, resp.\\
These graphs, due to Gray and Folkman, respectively, are edge-transitive but {\sc not} vertex-transitive. 
In fact, Tutte \cite[p. 67]{tutte} proved that any such graph {\it must} be bipartite.
Further, both the Gray and Folkman graphs are nonplanar.

In contrast, we know of many graphs, all vertex transitive, for which Conjecture BK
does hold, including cycles, 
hypercubes, 
complete bipartite graphs, degree-3 circulants, and the Heawood graph \cite{bk79, so-thesis, pck-circ}. Also, the Cartesian product of dispersable graphs is dispersable, and if $G$ is bipartite, then $mbt(G)$ is bounded above by the minimum sum of the page-maximum degrees over all book embeddings and trivially by $bt(G) \Delta(G)$; see \cite{ko-disp}.

In \cite{abgkp2018} Alam-et-al  proved Conjecture BK 
holds if $G$ is cubic planar  bipartite 3-connected and conjectured the condition, that the graph be 3-connected, could be relaxed. We confirm their conjecture by proving that all cubic planar bipartite {\it multigraphs} are dispersable, so Conjecture BK holds in this case. 

The paper is organized as follows.  Section 2 has some combinatorial preliminaries, including two ``folklore'' results for regular bipartite graphs---non-existence of a cutpoint and an ``entanglement'' lemma: any two cut-edges have vertices of attachment which are in different parts of the bipartition. Section 3 contains a proof of the main result.

\section{Combinatorial preliminaries}

Undefined graph terminology is standard \cite{harary}, \cite{diestel}. Let $|\cdot|$ denote cardinality.  

By a graph $G=(V,E)$ we mean a 1-dimensional simplicial complex; i.e., no loops and at most one edge between any pair of vertices; $G$ is a {\bf multigraph} if we allow {\bf parallel} edges, which are sets of more than one edge between a pair of vertices. Thus, every graph is a multigraph.  A {\bf digraph} is a multigraph whose edges have been {\bf oriented} (that is, assigned a direction).


A set of edges in a multigraph $G$ is a {\bf matching} if the edges are pairwise disjoint; a matching is {\bf perfect} if each vertex is in an edge of the matching.  Each page of a dispersable embedding of a regular graph is a perfect matching.


The following must be known but we haven't found it in the literature.

\begin{lemma}
Let $G$ be a bipartite, connected, and  $k$-regular multigraph, $k\geq3$.  Let $e', e''$ be two disjoint edges in $G$ such that $G \setminus \{e',e''\}$ contains a connected component $H$, and let $u', u''$ be the endpoints of $e', e''$ which are in $H$.  Then $u', u''$ are in distinct parts of the unique bipartition of $V(G)$.
\label{lm:entangle}
\end{lemma}
\begin{proof}
The {\it oriented degree} $\deg(v)$ of a vertex in a digraph is the number of edges in minus the number of edges out.
The sum of all oriented vertex degrees for any digraph is zero. 

Let $(H$, $u'$, $u'')$ be as in the statement; give $G$ and $H$ the edge-orientation induced from the bipartition by orienting the edges  from white to black, so  each vertex has oriented degree $\pm k$ according to whether it's black or white.
Put
$$W := V(H \setminus \{u',u''\})=V(H) \setminus \{u',u''\}.$$ We have
\begin{equation}
0 \;\;= \sum_{v \in V(H)} \deg_H(v) \,= \,\deg_H(u')+\deg_H(u'') + \sum_{v \in W} \deg_H(v).
\label{eq:par}
\end{equation}
For each  $v \in W$, $\deg_H(v) = \deg_G(v)$, so $\sum_{v \in W} \deg_H(v) \cong 0$ (mod $k$), hence, by (\ref{eq:par}), $\deg_H(u')+\deg_H(u'') \cong 0$. As $k \geq 3$, $2k-2 \ncong 0$ (mod $k$),
so the degrees of $u'$ and $u''$ must be of opposite sign, i.e., $u'$ and $u''$ are in different parts.
\end{proof}

By a similar argument, one can prove the following ``folklore'' result.

\begin{proposition}
Connected regular bipartite graphs have no cutpoints.
\label{pr:no-br}
\end{proposition}
\begin{proof}
Let $G$ be a connected $k$-regular bipartite graph. If $k=0$ or $1$, then $G$ is $K_1$ or $K_2$, resp., and there are no cutpoints. For $k=2$, bipartite or not, $G$ is a cycle and so has no cutpoint or bridge.  Suppose now $k \geq 3$. 
Let $u$ be a cutpoint of $G$ and let $\Gamma$ be a connected component of $G-u$.  

We define $J := \Gamma + u := G(V(\Gamma) \cup \{u\})$ to be the smallest subgraph of $G$ containing $\Gamma$ and $u$. Then again using oriented degrees,
\begin{equation}
0\;\; = \sum_{v \in V(J)}\deg_J(v)\, =\, \deg_J(u) + \sum_{v \in V(\Gamma)} \deg_G(v);
\end{equation}
the second summand is congruent to zero (mod $k$)
but $\deg_J(u) \ncong 0$ (mod $k$) so no cutpoint $u$ can exist.  
\end{proof}

Observe that for $k \geq 2$,
there is also no bridge.

\section{Cubic planar bipartite graphs}

Alam-et-al. proved the following result \cite[pp. 18--21]{abgkp2018}.

\begin{theorem}[\cite{abgkp2018}]
Let $G$ be a bipartite cubic planar $3$-connected graph.  Then $G$ has a subhamiltonian vertex ordering $\omega$ such that $mbt(G) = mbt(G,\omega) = 3$.
\label{th:alam}
\end{theorem}
\vspace{-.4 cm}
We say that such a $G$ is {\bf subhamiltonian dispersable}.
It is conjectured in \cite{abgkp2018} that $3$-connected is not needed, as we now show.
\begin{theorem}
Bipartite cubic planar multigraphs are subhamiltonian dispersable.
\label{th:bcpid}
\end{theorem}
\vspace{-.4 cm}
Our proof of Theorem \ref{th:bcpid} uses   multigraphs. 
The {\bf trivial cubic multigraph}, which we denote by $\Theta$, has 2 vertices joined by 3 parallel edges.

\begin{lemma}
A cubic multigraph $\neq \Theta$ with parallel edges cannot be $3$-connected.
\end{lemma}
\begin{proof}
If $G$ is cubic and $f_1,f_2$ are parallel edges joining vertices $a$ and $b$, then either (i) there is a third parallel edge $f_3$ and $G$ contains the trivial cubic graph as a connected component or (ii) the remaining two edges $f_3$ and $f_4$ at $a$ and $b$, resp., are a cutset for $G$, as they separate $\{a,b\}$ from the other vertices.
\end{proof}

\begin{center}
\includegraphics[width=0.7\linewidth]{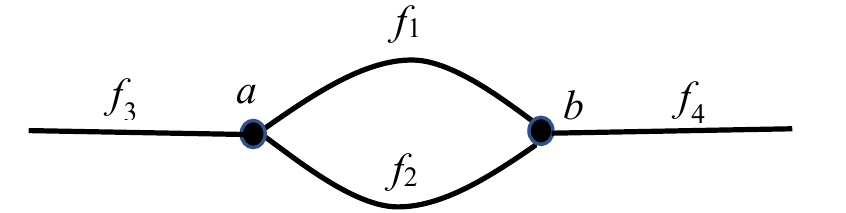}
\end{center}
\centerline{{\bf Figure 1} \ \ {\it Parallel edges in a cubic graph}}
\bigskip
\bigskip

Hence, {\it a nontrivial 3-connected cubic multigraph has no parallel edges and is a cubic} graph.
We are now ready for the proof of Theorem \ref{th:bcpid}.

\begin{proof}
Let $G$ be a bipartite cubic planar multigraph.  Assume without loss of generality (\cite{bk79}, \cite{so-thesis}) that $G$ is connected. 
We need to show that $G$ has a subhamiltonian vertex order $\omega$ and an edge-3-coloring $c$ such that in the outerplane drawing $(G, \omega)$, each color class of $c$ is crossing-free; that is, 
\begin{equation}
bt(G) = bt(G,\omega) \leq 2
\;\;\mbox{{\rm and}}\;\;mbt(G) = mbt(G, \omega) = 3.
\end{equation}

Suppose that $G$ is a counterexample with the minimum number of vertices.  Note that $G$ must be a nontrivial cubic multigraph as the required $\omega$ and $c$ do exist for $\Theta$, and $G$ cannot be $3$-connected  else it is dispersable by Theorem \ref{th:alam} which applies as there can be no parallel edges.  

Since $G$ is bipartite and cubic, by Proposition \ref{pr:no-br} it has no cutpoint and hence must be 2-connected.
As line-connectivity and point-connectivity coincide for cubic graphs \cite[p. 55, ex. 5.6 ]{harary}, there must be a cutset of two edges.  But the
two edges of such a cutset cannot share a vertex -- the third edge at a common vertex would be a bridge.
Hence, $G$ contains vertex-disjoint edges $e'=u'w',\; e''=u''w''$  (neither of which is a parallel edge) such that 
\begin{equation}
G \setminus \{e',e''\} = G_1 \cup G_2, \;\;u', u'' \in V(G_1), \;\;w', w'' \in V(G_2),
\end{equation}
where $G_1$ and $G_2$ are vertex-disjoint.  Further, by Lemma \ref{lm:entangle}, $u'$ and $u''$ are in distinct parts of the bipartition of $G_1$ induced by the bipartition of $G$, and similarly for $w', w''$ in the corresponding bipartition of $G_2$.

Let $H_1 = G_1 + e_1$, where $e_1 := u'u''$, $H_2 = G_2 + e_2$, where $e_2 := w'w''$.  Then $H_1$ and $H_2$ are bipartite cubic planar multigraphs, either of which can contain multiple edges even if $G$ has no parallel edges. As $|H_1|, |H_2| < |G|$, there exist cyclic orderings $\omega_1, \omega_2$ on $V(H_1)=V(G_1)$ and $V(H_2)=V(G_2)$, resp., and edge colorings $c_1, c_2$ which constitute dispersable subhamiltonian book embeddings of $H_1$ and $H_2$.  
Let $\lambda_1$ be  one of the two possible linear orderings obtained from $\omega_1$ which  starts at $u''$.  Similarly, let $\lambda_2$ be  one of the two possible linear orderings obtained from $\omega_2$ which  starts at $w''$. 

Define $\lambda$ a linear order on $V(G) = V(H_1) \cup V(H_2)$ by the equation $ \lambda := \lambda_1^{op} \,*\, \lambda_2$,
where $\lambda_1^{op}$ is the opposite (i.e., reverse) order to $\lambda_1$ and $*$ denotes concatenation.
Let $\omega$ be the cyclic order determined by $\lambda$.
Observe that $ u''$ and $w''$ are consecutive in the subhamiltonian ordering of $V(G)$ given by $\omega$.

By renaming the colors, we can assume that the colorings $c_1$ and $c_2$ both use colors $\{\alpha,\beta, \gamma\}$ and that $e_1$ and $e_2$ are both assigned $\gamma$.
Define $c$ an edge 3-coloring of $E(G)$ to be $c_1$ on $E(G_1)$, $c_2$ on $E(G_2)$, and put $c(e')=c(e'')=\gamma$. If $e \in E(G_1)$ is adjacent to either $e'$ or $e''$, then $e$ was adjacent to $e_1$, so $c(e) \neq c(e')$ and similarly for $G_2$. Observe that $e''$ has no crossings in $(G,\omega)$. 

\begin{center}
\includegraphics[
width=0.43\linewidth]{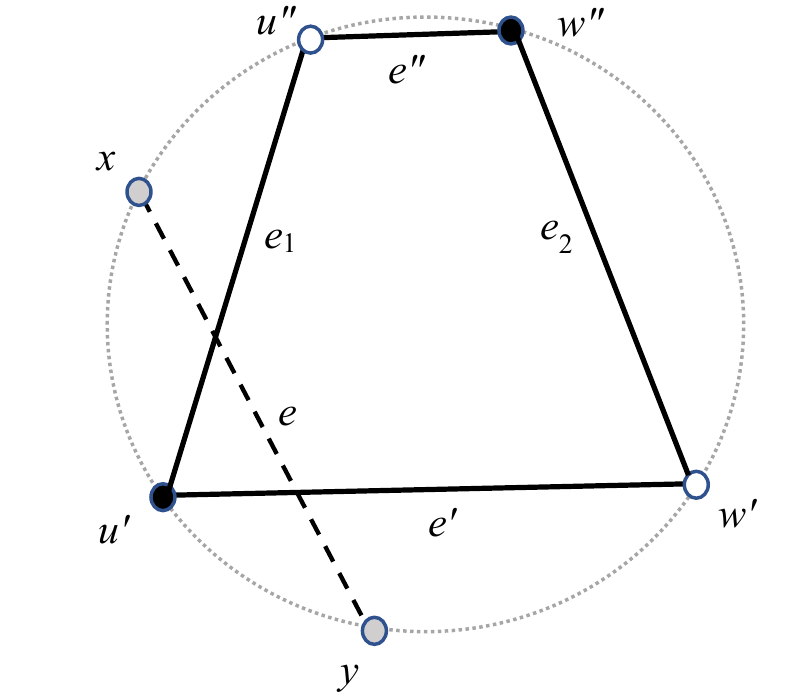}
\end{center}
\centerline{{\bf Figure 2} \ \ {\it Subhamiltonian dispersable embedding of $G$}}
\bigskip
\bigskip

We claim $e'$ crosses no edge colored $\gamma$.  For if $e=xy$ crosses $e'$ in the drawing $(G,\omega)$, then $x$ and $y$ both belong to $G_1$ or both belong to $G_2$; else $\{e', e''\}$ wouldn't be a cutset.
Thus $xy$ must cross $e_1$ or $e_2$, hence $c(xy) \neq \gamma$. So $(\omega,c)$ is a dispersable subhamiltonian book embedding of $G$ which is impossible.
\end{proof}

\section{Discussion}


Since the examples of non-dispersable bipartite graphs given in Alam-et-al. are nonplanar and not vertex transitive, 
natural questions arise. They ask \cite[question 6]{abgkp2018}: Does the BK Conjecture in \cite{bk79} 
hold for planar bipartite graphs? 

All examples we know support the following conjecture which would include a positive answer to \cite[questions 2 and 4]{abgkp2018}:
{\it Every vertex-transitive graph $G$ with $\Delta(G) = r$ has $mbt(G) = r$ {\rm or} $\,r{+}1$ according to whether or not $G$  is bipartite.}  However, some condition on $G$ is needed as Barat et al.\cite{bmw2006} show that if $\Delta(G) \geq 9$, bounded degree graphs cannot have bounded book thickness.

\Addresses
\bigskip

\subsection*{Postscript}
The previous paper [2] proves the conjecture of Alam et al. [1] made in their 2018 conference paper which we discovered in 2019.  We worked during the Spring of 2020 and on July 20, 2020, our paper was submitted for publication in the form above.  After 8 months, on March 26, 2021, we contacted the assistant to the managing editor of the journal and were told that no referee had been found.  We then supplied names, affiliations, and emails for several (we thought obvious) choices for possible referees, which were accepted with thanks.  On May 30, 2021, the journal rejected the paper. The referee's report said that the problem had been solved and included a link to the now-published expanded version [3] of Alam et al.'s results.

Indeed, on March 3, 2020, Alam, et al. first submitted a version of the conference paper [1] to the Elsevier journal, Theoretical Computer Science and submitted the second version [3] on Jan. 23, 2021, and it appeared electronically on Feb. 4, 2021. This new journal version [3], with a sixth coauthor, Dujmovi\'{c},  includes  a proof of the  conjecture in [1].
  
Another proof of the Alam et al. conjecture was recently given by Shao, Liu, and Li [4].  The paper [4] uses a different method, while [2] and the portion of [3], which proves the conjecture, use very similar techniques.

We hope the comparison of these approaches will be helpful.
\\
\begin{flushleft}
{\large {\bf References}}
\bigskip

[1]	J. Md. Alam, M. A. Bekos, M. Gronemann, M. Kaufmann \& S. Pupyrev, On dispersable book embeddings, Graph-theoretic Concepts in Computer Science (44th International Workshop, WG 2018, Cottbus, Germany),  A. Brandest\"{a}dt, E. K\"{o}hler, \& K. Meer, Eds., LNCS 11159 (2018) 1–14, Springer, Cham, Switzerland.\\
\vspace{0.2cm}
[2]	P. C. Kainen \& S. Overbay, Cubic planar bipartite graphs are dispersable, July 20, 2020 (7 pages).\\
\vspace{0.2cm}
[3]	 J. Md. Alam, M. A. Bekos, V. Dujmovi\'{c}, M. Gronemann, M. Kaufmann \& S. Pupyrev, On dispersable book embeddings, Theoretical Computer Science, 861 (2021), 1--22.\\
\vspace{0.2cm}
[4]	Z. Shao, Y. Liu \& Z. Li, Bipartite cubic planar graphs are dispersable, arXiv:2107.00907v1, July 2, 2021, (10 pages).
\end{flushleft}
\end{document}